\documentstyle[12pt]{article}
\begin{document}
\topmargin-.75cm
\baselineskip16.125pt
\parskip.5\baselineskip

\renewcommand{\thesection}{\ifnum\value{section}=0 0\else
\arabic{section}\fi}

\newtheorem{definition}{Definition $\!\!$}[section]
\newtheorem{prop}[definition]{Proposition $\!\!$}
\newtheorem{lem}[definition]{Lemma $\!\!$}
\newtheorem{corollary}[definition]{Corollary $\!\!$}
\newtheorem{theorem}[definition]{Theorem $\!\!$}
\newtheorem{example}[definition]{Example $\!\!$}
\newtheorem{remark}[definition]{Remark $\!\!$}

\newcommand{\nc}[2]{\newcommand{#1}{#2}}
\newcommand{\rnc}[2]{\renewcommand{#1}{#2}}
\nc{\bpr}{\begin{prop}}
\nc{\bth}{\begin{theorem}}
\nc{\ble}{\begin{lem}}
\nc{\bco}{\begin{corollary}}
\nc{\bre}{\begin{remark}}
\nc{\bex}{\begin{example}}
\nc{\bde}{\begin{definition}}
\nc{\ede}{\end{definition}}
\nc{\epr}{\end{prop}}
\nc{\ethe}{\end{theorem}}
\nc{\ele}{\end{lem}}
\nc{\eco}{\end{corollary}}
\nc{\ere}{\end{remark}}
\nc{\eex}{\end{example}}
\nc{\epf}{\hfill\mbox{$\Box$}}
\nc{\ot}{\otimes}
\nc{\LAblp}{\mbox{\LARGE\boldmath$($}}
\nc{\LAbrp}{\mbox{\LARGE\boldmath$)$}}
\nc{\blp}{\mbox{\boldmath$($}}
\nc{\brp}{\mbox{\boldmath$)$}}
\nc{\LAlp}{\mbox{\LARGE $($}}
\nc{\LArp}{\mbox{\LARGE $)$}}
\nc{\Llp}{\mbox{\Large $($}}
\nc{\Lrp}{\mbox{\Large $)$}}
\nc{\llp}{\mbox{\large $($}}
\nc{\lrp}{\mbox{\large $)$}}
\nc{\lbc}{\mbox{\Large\boldmath$,$}}
\nc{\lc}{\mbox{\Large$,$}}
\nc{\Lall}{\mbox{\Large$\forall$}}
\nc{\bc}{\mbox{\boldmath$,$}}

\nc{\Section}{\setcounter{definition}{0}\section}
\renewcommand{\theequation}{\thesection.\arabic{equation}}

\newcounter{c}
\renewcommand{\[}{\setcounter{c}{1}$$}
\newcommand{\etyk}[1]{\vspace{-7.4mm}$$\begin{equation}\label{#1}
\addtocounter{c}{1}}
\renewcommand{\]}{\ifnum \value{c}=1 $$\else \end{equation}\fi}

\newcommand{\dowod}{\noindent{\bf Proof:} }
\newcommand{\Sp}{{\rm Sp}\,}
\newcommand{\Mor}{\mbox{$\rm Mor$}}
\newcommand{\skrA}{{\cal A}}
\newcommand{\Phase}{\mbox{$\rm Phase\,$}}
\newcommand{\id}{{\rm id}}
\newcommand{\diag}{{\rm diag}}
\newcommand{\inv}{{\rm inv}}
\newcommand{\ad}{{\rm ad}}
\newcommand{\poi}{{\rm pt}}
\newcommand{\Dim}{{\rm dim}\,}
\newcommand{\Ker}{{\rm ker}\,}
\newcommand{\Mat}{{\rm Mat}\,}
\newcommand{\Rep}{{\rm Rep}\,}
\newcommand{\Fun}{{\rm Fun}\,}
\newcommand{\Tr}{{\rm Tr}\,}
\newcommand{\supp}{\mbox{$\rm supp$}}
\newcommand{\half}{\frac{1}{2}}

\newcommand{\skrF}{{\cal F}}
\newcommand{\skrD}{{\cal D}}
\newcommand{\skrC}{{\cal C}}

\newcommand{\ttimes}{\mbox{$\hspace{.5mm}\bigcirc\hspace{-4.9mm}
\perp\hspace{1mm}$}}
\newcommand{\Ttimes}{\mbox{$\hspace{.5mm}\bigcirc\hspace{-3.7mm}
\raisebox{-.7mm}{$\top$}\hspace{1mm}$}}
\newcommand{\Cstar}{$^*$-}

\newcommand{\bbr}{{\bf R}}
\newcommand{\bbz}{{\bf Z}}
\newcommand{\Ci}{C_{\infty}}
\newcommand{\Cb}{C_{b}}
\newcommand{\fa}{\forall}
\newcommand{\te}{\exists}
\newcommand{\rrr}{right regular representation}
\newcommand{\wrt}{with respect to}

\newcommand{\qg}{quantum group}
\newcommand{\qgs}{quantum groups}
\newcommand{\cs}{classical space}
\newcommand{\qs}{quantum space}
\newcommand{\po}{Pontryagin}
\newcommand{\ch}{character}
\newcommand{\chs}{characters}
\newcommand{\hg}{\hat{G}}
\newcommand{\be}{\begin{equation}}
\newcommand{\ee}{\end{equation}}
\newcommand{\al}{\alpha}
\newcommand{\bet}{\beta}
\newcommand{\Lam}{$\Lambda$}
\newcommand{\La}{\Lambda}
\newcommand{\lam}{$\lambda$}
\newcommand{\la}{\lambda}
\newcommand{\ep}{\epsilon}
\newcommand{\Om}{\Omega}
\newcommand{\lr}{\longleftrightarrow}
\newcommand{\g}{\Gamma}
\newcommand{\gt}{{\widetilde{\Gamma}}}
\newcommand{\D}{\Delta}
\newcommand{\si}{\sigma}
\newcommand{\f}{\varphi}
\newcommand{\eps}{\varepsilon}
\newcommand{\w}{\omega}
\newcommand{\C}{{\bf C }}
\newcommand{\R}{{\bf R }}
\newcommand{\bea}{\begin{eqnarray}}
\newcommand{\eea}{\end{eqnarray}}
\newcommand{\ba}{\begin{array}}
\newcommand{\ea}{\end{array}}
\newcommand{\ber}{\begin{eqnarray}}
\newcommand{\ear}{\end{eqnarray}}
\newcommand{\U}{\Upsilon}
\newcommand{\der}{\mbox{d}}
\newcommand{\de}{\mbox{{\rm det}}}

\def\inbar{\,\vrule height1.5ex width.4pt depth0pt}
\def\IC{\relax\,\hbox{$\inbar\kern-.3em{\rm C}$}}
\def\otc{\otimes_{\IC}}
\def\ra{\rightarrow}
\def\ota{\otimes_ A}
\def\otza{\otimes_{ Z(A)}}
\def\otc{\otimes_{\IC}}
\def\h{\rho}
\def\x{\zeta}
\def\th{\theta}
\def\s{\sigma}
\def\t{\tau}
\def\st{\stackrel}
\def\ot{\otimes}
\def\dq{\Delta_q}
\title{{\bf 
ON THE DRINFELD TWIST\\
FOR QUANTUM $sl(2)$ }}
\author{Ludwik D\c{a}browski\\
\normalsize SISSA, Via Beirut 2-4, Trieste, Italy.\\
\normalsize \sc e-mail: dabrow@sissa.it\vspace*{5mm} \\
Fabrizio Nesti \\ 
\normalsize SISSA, Via Beirut 2-4, Trieste, Italy.\\
\normalsize \sc e-mail: nesti@sissa.it\vspace*{5mm}\\
Pasquale Siniscalco\\
\normalsize SISSA, Via Beirut 2-4, Trieste, Italy.\\
\normalsize\sc e-mail: sinis@sissa.it \vspace*{5mm}}
\date{ }
 \maketitle
\thispagestyle{empty}
\setcounter{page}{0}
\begin{abstract}
\noindent
An isomorphism, up to a twist, between the quasitriangular quantum 
{\mbox enveloping} {\mbox algebra} $U_h(sl(2))$ and the (classical)
$U(sl(2))[[h]]$ ~is discussed. 
The universal {\mbox twisting} element $\cal F$ is given up to the 
second order in the deformation parameter $h$.
\end{abstract}
\vfill
\centerline{SISSA 130\slash 96\slash FM}
\newpage\parindent5mm\baselineskip16.5pt\parskip.65\baselineskip

\Section{Introduction}
In 1989 Drinfeld showed by cohomological arguments that, 
as a formal series in a deformation parameter $h$,  
all the quantum symmetries (quasitriangular Hopf algebras) 
$U_h(g)$, where $g$ is a semisimple Lie algebra, 
are isomorphic with $U(g)[[h]]$, up to a twist $\cal F$ \cite{dr, dri}.
He also posed a problem (cf.\cite{dri}) to find a concrete pair $(m, \cal F)$,
consisting of an isomorphism $m$ and a universal twisting element $\cal F$.
This turns out to be a formidable task, which as far as we know, is not yet solved
in general.
The only case when it has been performed concerns the q-deformed 
Heisenberg algebra ${\cal H}_q(1)$ \cite{bgst}.
The next important case to be investigated 
is the quantum deformation of $sl(2)$
(as a matter of fact ${\cal H}_q(1)$ can be obtained from it by a contraction).
As far as $U_h(sl(2))$ is regarded, a candidate for the isomorphism $m$ is actually 
known \cite{cz}. 
Also, a series of related particular matrix solutions for the twist element $\cal F$
were reported, 
namely $\cal F$ in the representations $\frac{1}{2} \ot j $, where 
$\frac{1}{2} $ denotes the fundamental representation and 
$j$ denotes the irreducible $(2j+1)$-dimensional  representation of 
$sl(2)$ \cite{c,z}, (see also \cite{enge}).
Moreover, in \cite{cgz} a sort of a `semi-universal' form of
$\cal F$ has been given,
i.e. the expression for $(\frac{1}{2} \ot \id ) (\cal F)$.
However, the universal element $\cal F$ itself has not been known 
beyond the first order in the deformation parameter $h$
(the first order coefficient being given by the classical r-matrix $r$).
In this letter, we investigate and report the solution up to the second order 
in $h$. In the subsequent sections we separately discuss the problem on the levels 
of algebra, Hopf algebra and quasitriangular Hopf algebra.

It is worth to mention that evaluating $\cal F$ 
in the representation $\rho_L\otimes \rho_L$,
where $\rho_L$ is the representation of $sl(2)$ in terms of the
left-invariant vector fields on $SL(2)$,
one obtains a quantization of the Lie-Poisson bracket on $SL(2)$ given 
by $r$ \cite{tak}. 
In particular, the second coefficient of $(\rho_L\otimes \rho_L ) (\cal F)$ 
provides an interesting second order (bi)differential operator on $SL(2)$.

\parindent0pt\baselineskip16.125pt\parskip.5\baselineskip
\section{Algebra level}
 We start by specifying our conventions about Lie algebra $sl(2)$.
The generators are $H, E, F$ with the commutation relations: 
\be 
\label{comm}
[H,E]=E,~~~~ [H,F]=-F,~~~~ [E, F]=H~. 
\ee
As a consequence we have the following exchange relations between 
any polynomial $\phi(H)$ in $H$ 
and the 
powers of $E$ and $F$:
\bea \phi(H)E^n & = & E^n \phi(H+n)~, \nonumber \\
      \phi(H)F^n & = & F^n \phi(H-n) ~.
\eea
The quadratic Casimir element in the universal enveloping algebra $U(sl(2))$ is
\be I=2EF+ H(H-1) =2FE +H(H+1)\doteq j(j+1)~. \label{cas} \ee 
A possible basis for the enveloping algebra is provided by the set 
$\{H^l E^m F^n\}$, 
but using the relations (\ref{cas}) we can pass to the basis given by 
$\{H^a I^b E^c \oplus H^r I^s F^t\}$. 
This basis will be more suitable for our computations.\\
Next, the generators $ J^{+}, J^-, J^0$ of the q-deformed 
algebra obey the following commutation relations:
 \be 
[J^0,J^+]=J^+,~~~~ [J^0,J^-]=-J^-,~~~~ [J^+,J^-]=\frac{1}{2}[2J^0]~. 
 \label{qcomm} \ee
where $[x]$, the {\em q-analogue} of $x$, is defined as:
\be 
[x]= \frac{q^x-q^{-x}}{q-q^{-1}} ~.
\ee 
The `deforming maps' introduced in \cite{cz}, provide (cf.\cite{d}) an isomorphism 
$m$ between $U_{h}(sl(2))$ and $U(sl(2))[[h]]$
which is given by mapping the generators 
$J^0, J^+, J^-$ to the following combinations of $H, E, F$ 
\be 
\label{iso}
J^0 \ra H, ~~~J^{+} \ra \phi^+ E,~~~ J^- \ra \phi^- F =F \phi^+ ~,
\ee 
 where 
\be 
\label{iso1}
\phi^{\pm} = \sqrt{\frac{[j\pm H][1+j\mp H]}{(j\pm H)(1+j \mp H )} } ~.
\ee
We remark that (\ref{iso1}) is a well defined expression, 
as the inverse and square root operations are admissible in the $h$-adic topology.
In fact, with $q=e^{h}$,
we can write the expansion in $h$ up to the second order as 
 \be
\label{phie}
\phi^{\pm}  = 
 1 + \frac{1}{12} h^2 \left(2I +2H(H\mp 1) -1\right) + o(h^3)
\doteq 1+ h^2 \phi^{\pm}_2 + o(h^3) ~.
 \ee
It will be useful to mention \cite{Kas}, that any other isomorphism $m'$ 
can differ 
at most by a similarity via an invertible element $M\in U(sl(2))[[h]]$, 
i.e. \be
\label{sim}
m' = M m M^{-1} ~.
\ee
We conclude this section with few remarks.
Note that (\ref{iso}) is in fact valid also for the $^*$-algebras
$U_{h}(su(2))$ and $U(su(2))[[h]]$, since it fulfills the relevant hermicity condition.
In this respect, more general isomorphisms belonging to the
one-parameter family
introduced in \cite{cz} do not satisfy such a hermicity requirement.
In addition, they are not suitable for our purposes since the coefficients 
of the expansion in $h$ are not polynomial in the generators.
\section{Hopf algebra level}  
The enveloping algebra $U(sl(2))[[h]]$ with relations (\ref{comm})
when equipped with the  usual coproduct
\be \D(x) = 1\otimes x + x \otimes 1 ~, ~~~\forall x\in sl(2) ~,
\ee
becomes a Hopf algebra. In the quantum case, the coproduct in $U_{h}(sl(2))$ 
is defined as:
\bea \Delta_q(J^0) &= &1\otimes J^0 + J^0 \otimes 1 ~, 
     \nonumber \\
      \dq(J^{\pm}) & = &J^{\pm}\otimes q^{J^0}+q^{-J^0} \otimes J^{\pm}~. 
\eea
(The counit and coinverse are not needed for our purposes).\\
The main part of Drinfeld Theorem guarantees that these two classical and quantum 
coproducts are related via a twist by an invertible 
${\cal F} \in \left( U(sl(2))\otimes U(sl(2))\right) [[h]]$.\\
More precisely, defining
\be  \widetilde{\dq} \doteq (m\ot m) \circ \dq \circ m^{-1} \ee 
we have
\be \widetilde{\dq}(x) = {\cal F} \D (x){\cal F}^{-1},~~~\forall x
 \in U(sl(2))[[h]] \ee 
 It is sufficient (and necessary)
to verify this equation by substituting for $x$ the image by $m$ of the generators  
 $J^0,J^+,J^-$. \\ 
We remark that there is no loss of generality in restricting ourselves 
to a specific isomorphism (\ref{iso}). 
Indeed,  had we used another isomorphism $m'$, it turns out from (\ref{sim}) 
that the corresponding $\cal{F}'$ would be given by 
$(M\ot M)\widetilde{\dq}(M){\cal F}$.

As it is known (cf. \cite{tak}) a particular solution up to first order in 
$h$ is just ${\cal F}=1 + h r$, where 
\be 
\label{crm}
r = F\ot E -E \ot F \ee
is the standard classical r-matrix. More generally and up to order two in $h$ we write
\be 
{\cal F} ={\cal F}_0 + h {\cal F}_1 + h^2 {\cal F}_2 + o(h^3)~,
\ee  
with ${\cal F}_i$ belonging  to $U(sl(2)) \ot U(sl(2))$.\\ 
Using (\ref{phie}), we obtain the following coupled system of equations to solve 
by recursion:
\bea  
&& [{\cal F}_0,\D H] =0 ~, \nonumber \\ 
 && [{\cal F}_0,\D E] =0 ~, \nonumber \\ 
 &&  [{\cal F}_0,\D F] =0 ~. \label{F0} \eea
\bea
 && [{\cal F}_1,\D H] =0 ~,\nonumber \\ 
 && [{\cal F}_1,\D E] =(E\ot H - H \ot E){\cal F}_0 ~,\nonumber \\ 
 &&  [{\cal F}_1,\D F] = (F\ot H - H \ot F) {\cal F}_0 ~. \label{F1} \eea   
\bea
 && [{\cal F}_2,\D H] =0 ~,\nonumber \\
 &&  [{\cal F}_2, \D E  ] = (E\ot H - H \ot E) {\cal F}_1 
 - {\cal F}_0\D \phi^+_2 \D E	\nonumber \\ 
 &&~~~~~~~~~~~~~ + ~( \frac{1}{2} E\ot H^2+ 
\frac{1}{2} H^2\ot
 \phi^+_2 E + \phi^+_2 E \ot 1 + 1\ot \phi^+_2 E){\cal F}_0  ~, \nonumber \\
 && [{\cal F}_2, \D F  ] = (F\ot H - H \ot F) {\cal F}_1 
	-{\cal F}_0\D \phi^-_2 \D F  \nonumber \\
 && ~~~~~~~~~~~~ + ~( \frac{1}{2} F\ot H^2+ \frac{1}{2} H^2\ot
 \phi^-_2 F + \phi^-_2 F \ot 1 + 1\ot \phi^-_2 F){\cal F}_0  ~. \label{F2} \eea
Besides ${\cal F}_0=1\ot 1$, any 
arbitrary polynomial $f_0$ in the variables $(I\ot 1, 1\ot I, \D I)$ 
satisfies equations (\ref{F0}).
Due to linearity of the equation we can write then: 
\be {\cal F}_0=1\ot 1+f_0~.\ee
As regards ${\cal F}_1$, besides the solution   
$\widetilde{{\cal F}_1} = r$ of the equations (\ref{F1}) (with $f_0=0$), 
a solution for the general case is given by 
\be {\cal F}_1=\widetilde{{\cal F}_1}(1\ot 1+f_0) +f_1~,
\ee 
with $f_1$ being a solution of (\ref{F0}).\\
Similarly for ${\cal F}_2$: if one finds a particular solution 
$\widetilde{{\cal F}_2}$ of 
(\ref{F2}) (with $f_0=f_1=0$), the most general one is given by 
\be {\cal F}_2=\widetilde{{\cal F}_2}(1\ot 1+f_0) + 
\widetilde{ {\cal F}_1}f_1 + f_2~,\ee with $f_2$ 
solution of (\ref{F0}).\\
The possibility of adding {\em pure kernel} 
(i.e. satisfying the homogeneous equations (\ref{F0}))
 terms $f_1$ and $f_2$
comes from the fact that the last two equations for ${\cal F}_1$ and 
${\cal F}_2$ are 
linear non homogeneous, whose associated homogeneous ones are 
the last two equations in (\ref{F0}).\\

Now we proceed to exhibit the aforementioned particular solutions 
$\widetilde{{\cal F}_i}$ of this set of equations.
In $U(sl(2)) \ot U(sl(2))$ we use the basis 
$$\{H^{a_1}I^{b_1}E^{c_1}
\oplus H^{r_1}I^{s_1}F^{t_1}\} \ot \{H^{a_2}I^{b_2}E^{c_2}
\oplus H^{r_2}I^{s_2}F^{t_2}\}~.$$ 
In order to simplify the notation, for any $x \in U(sl(2))$ we set 
$x_1= x\ot 1$, $x_2=1\ot x$. From $[{\cal F}_i,\D H] =0$,
for all $i$, it is easily seen that any 
$ {\cal F}_i$ is of the form 
${\cal F}_i =a_{il}E_1^l F_2^l  + b_{il} F_1^l E_2^l $,
where $a_{il}$ and $b_{il}$ are polynomials in $H_1, H_2, I_1, I_2$.\\
We've already mentioned that $\widetilde{{\cal F}_0}=1$ is a solution for 
equations (\ref{F0}).\\
Next we pass to the first order term. For simplicity we
drop the index $i=1$ in the following formulae and define 
\bea && \delta_1(a_k)= 
a_{k}(H_1,H_2,I_1,I_2)-a_{k}(H_1\! -\! 1,H_2,I_1,I_2)~, \nonumber \\
&& \delta_2(a_k)= 
a_{k}(H_1,H_2,I_1,I_2)-a_{k}(H_1,H_2\! -\! 1,I_1,I_2)~,\nonumber \eea 
and similarly for $b_k$.
The equations (\ref{F1}) give the following system of coupled partial difference 
equations for the coefficients $a_l$ and $b_l$:
\bea 
&& \delta_1(a_{n-1})=-\frac{1}{2}(I_2+H_2-H_2^2)\delta_2(a_n) 
+(nH_2 + \frac{n^2\! -\! n}{2}) a_n ,\nonumber \\
&& \delta_2(a_{n-1})=-\frac{1}{2}(I_1-H_1-H_1^2)
\delta_1(a_{n}(H_1\! +\! 1,H_2\! -\! 1)) +
(nH_1 - \frac{n^2\! -\! n}{2}) a_{n}(H_1,H_2\! -\! 1), \nonumber \\ 
&& \delta_1(b_{n-1})=-\frac{1}{2}(I_2-H_2-H_2^2)
\delta_2(b_{n}(H_1\! -\! 1,H_2\! +\! 1)) +
(nH_2 -\frac{n^2\! -\! n}{2}) b_{n}(H_1\! -\! 1,H_2), \nonumber \\
&& \delta_2(b_{n-1})=-\frac{1}{2}(I_1+H_1-H_1^2)\delta_1(b_{n}) 
+(nH_1 +\frac{n^2\! -\! n}{2}) b_{n}, \eea
for any $n\ge 2$, whereas for $n=1$ we have:
\bea 
&& \delta_1(a_0 + b_0)=-\frac{1}{2}(I_2+H_2-H_2^2)
\delta_2(a_1) + H_2 a_1 +H_2 ~,\nonumber \\
&& \delta_1(a_0+b_0)=-\frac{1}{2}(I_2-H_2-H_2^2)
\delta_2(b_1(H_1\! -\! 1,H_2\! +\! 1)) +H_2  b_{1}(H_1\! -\! 1,H_2)-H_2~, \nonumber \\
&& \delta_2(a_0+b_0)=-\frac{1}{2}(I_1-H_1-H_1^2)
\delta_1(a_{1}(H_1\! +\! 1,H_2\! -\! 1)) +H_1 a_{1}(H_1,H_2\! -\! 1)+H_1 ~,\nonumber \\
&& \delta_2(a_0+b_0)=-\frac{1}{2}(I_1+H_1-H_1^2)\delta_1(b_{1}) 
+H_1 b_{1}-H_1~. \eea 
In order to find a particular solution of this system of equations, 
one can fix a couple 
$\{N,K\}$ such that $a_n=b_k=0$, $\forall n\ge N$ and $\forall k \ge K$, in 
order to set the maximum degree 
for the polynomials in $E_1^l F_2^l$ and $E_2^l F_1^l$, and then solve 
recursively the equations 
for the lower degree terms by partial finite integration.\\
By making a {\em minimal} choice, putting $a_n=b_n=0$, for any 
$n\ge 2$, we recover the solution:
\be \widetilde{{\cal F}_1}= r~,
\ee
with $r$ given by (\ref{crm}).
Consistently with what we explained in the previous section, 
had we decided to fix our cut-off at 
higher degree terms we would have adjoined to $\widetilde{{\cal F}_1}$ some 
$f_1$ solution of the {\em pure kernel} part.\\
As regards ${\cal F}_2$, the structure of the equations for $a_l$ and $b_l$ 
remains unchanged for $n\ge 3$, whereas for $n=\{2,1\}$ some extra term 
appear, due to $\phi^+_2$ and $\phi^-_2$.\\
We skip the explicit (and lengthy) form of them, and we just give the 
expression for a particular solution:
\bea \widetilde{{\cal F}_2} & = & \frac{1}{2} (I\ot H^2+H^2\ot I)+
\frac{1}{3}(E\ot HF-HE\ot F+HF\ot E-F\ot HE) 
 \nonumber \\
 &  & +\frac{1}{6}H\ot H(1-3P) -\frac{11}{24}P +\frac{1}{2}\left( 
(1+P)^2-1-2I\ot I \right)~, \eea
where 
\be
\label{P}
P=2(E\ot F+F\ot E +H\ot H) 
\ee
is the Cartan-Killing metric.\\
Applying representations of $sl(2)$ we can obtain explicit matrix
expressions for $\widetilde{{\cal F}}$. 
It turns out that our particular solution $\widetilde{{\cal F}}$,
when composed with $\frac{1}{2} \ot \id $,
reproduces the semi-universal solution presented in \cite{cgz} 
in terms of $2\times 2$ matrices with coefficients in $U(sl(2))[[h]]$
(up to the second order in $h$).
Thus, as a consequence it also coincides with the matrix solutions 
in the representations $\frac{1}{2} \ot j $.\\
We remark that in the literature one may find often other properties of the twisting
element ${\cal F}$. For instance, ${\cal F}$ may be supposed to satisfy the \\
i) `normalization' condition 
\be 
(\varepsilon \ot id)({\cal F}) =
(id \ot \varepsilon )({\cal F}) = 1~, 
\ee 
sometimes also expressed as 
${\cal F}(x,0) = {\cal F}(0,y) =1$. With the standard definition of 
counit $\varepsilon$ this implies ${\cal F}_0=1$, i.e. $f_0=0$.\\
ii) unitarity condition $\sigma({\cal F}){\cal F}=1$.
In our case ${\cal F}$ fulfills this condition in the particular 
representation $\frac{1}{2}\ot \frac{1}{2}$, but not in general.\\ 
iii) condition $( {\cal F}\ot \id )(\D\ot\id ){\cal F}  = (\id \ot{\cal F} )
(\id \ot\D){\cal F}$. In our case $ {\cal F}$ does not fulfill it, not 
even in a representation (except the trivial one). We remark that this 
condition is a stronger requirement with respect to the coassociativity of 
the twisted coproduct, which in our case follows directly from the definition. \\

\section{Quasitriangular Hopf algebra level}
 From the Drinfeld theorem, the quantum universal R-matrix 
\be
R_q = q^{2J^0\ot J^0} \sum_{n=0}^{\infty} 
q^{-n(n-1)/2}\frac{2^n(1-q^{-2})^n} {[n]!} (q^{J^0}J^+ \ot q^{-J^0}J^-)^n ~,
\ee
and the undeformed universal R-matrix, though not the simple $1\ot 1$ but rather, 
\be 
R=q^P~, 
\ee
with $P$ given by (\ref{P}),
should be related by the isomorphism up to a twist.
Thus, setting 
\be
\widetilde{R_q} \doteq (m\ot m)(R_q) ~,
\nonumber 
\ee
 ${\cal F}$ is supposed to verify the equation 
\be 
\widetilde{R_q}  {\cal F}=\sigma({\cal F})R \label{req} ~,
\ee
where $\sigma$ is the flip operator and $[n]!\doteq [n][n-1]\dots 1$.\\
We have the following expansions
\bea
R &=& 1+ h R^{(1)} +h^2 R^{(2)}+ o(h^3)\nonumber \\
&=& 1+h (2 E\ot F+2F\ot E + 2H \ot H) \nonumber	\\
& & + h^2 (-H\ot H - 2 E \ot F - 2 F\ot E + 2 E^2\ot F^2 +2 F^2\ot E^2 
 \nonumber \\
& & -2 E \ot HF -2 HF\ot E + 2 F\ot HE +2 HE\ot F 
 +4 HE\ot HF +4 HF\ot HE  \nonumber \\
& & +3 H^2\ot H^2 + I\ot I  - I\ot H^2 -H^2\ot I )~,	
 \eea
\bea
\widetilde{R_q} &=& 1+ h R^{(1)}_q +h^2 R^{(2)}_q+ o(h^3)\nonumber \\
&=& 1+h (4 E\ot F+ 2H \ot H)\\
& & + h^2 (2 H^2\ot H^2 -4 E \ot F -4 E \ot HF + 8 E^2\ot F^2 +4 HE\ot F 
+ 8 HE\ot HF)~.	\nonumber
\eea
At the zero-order in $h$, choosing $f_0 =0$, (\ref{req}) is identically 
satisfied ($1=1$).\\ At the order one we have the following equation: 
\be \sigma(f_1) -f_{1} = R_q^{(1)} -R^{(1)} -\left(\sigma(\widetilde{{\cal F}_1}) -
\widetilde{{\cal F}_1}\right)~. \label{r1}\ee
It comes from direct computations that the right-hand-side is zero, 
which implies that
$f_1$ must be symmetric.\\
At the second order we obtain
\bea \sigma(f_2) -f_{2} 
  & = & \! \widetilde{{\cal F}_2}-\sigma 
 (\widetilde{{\cal F}_2}) -\sigma( \widetilde{{\cal F}_1}) R^{(1)} 
 +R_q^{(1)} \widetilde{{\cal F}_1} +R_q^{(2)}- R^{(2)}~. \label{r2} \eea 
Again the right-hand-side is zero, and hence also $f_2$ must be symmetric.\\
Since, in particular, $f_1$ and $f_2$ can be equal to zero, we have that 
our particular solution $\widetilde{{\cal F} }$ satisfies (\ref{req}).
\section{Conclusions}
In accordance with the theorem of Drinfeld, 
we have exhibited an isomorphism from $U_{h}(sl(2))$ to $U(sl(2))[[h]]$
and (up to the second order in $h$)
a class of universal twisting elements ${\cal F} \in \left( U(sl(2)) \ot 
U(sl(2)) \right) [[h]]$.
Such $\cal F$ perform a gauge transformation (twist) from the ordinary coproduct 
and from the universal R-matrix $R=q^P$ in $U(sl(2))[[h]]$ 
to their quantum counterparts in $U_{h}(sl(2))$.

We have identified a particular universal element
$\widetilde{{\cal F}}$ in this class which, 
after applying the representation $\frac{1}{2}$ to its first leg,
coincides with the `semi-universal' solution in \cite{cgz} 
(up to the second order in $h$). 
Consequently, it also coincides with the known matrix solutions 
in the representations $\frac{1}{2} \ot j $.

The computation of the higher order terms, with the help of `Mathematica',
is in progress.

\end{document}